\documentclass[12pt]{amsart}
\usepackage{amssymb}

\newtheorem{theorem}{Theorem}[section]

\newtheorem{lemma}{Lemma}[section]

\title[Binary recursive integer sequences]{Integers represented by binary recursive sequences}

\date{\today}

\author{Lajos Hajdu}
\address{Institute of Mathematics, University of Debrecen,\newline
\indent P. O. Box 400, H-4002 Debrecen, Hungary \newline
\indent and HUN-REN DE Equations, Functions, Curves and their Applications Research Group}
\email{hajdul@science.unideb.hu}

\author{Rob Tijdeman}
\address{Mathematical Institute, Leiden University, The Netherlands}
\email{tijdeman@ziggo.nl}

\thanks{Research supported in part by the HUN-REN Research Network and by the NKFIH grant 130909.}

\begin{document}

\begin{abstract}
This paper is the continuation of \cite{htl}, where we deal with Lucas sequences. Here we study integers represented by integer sequences which satisfy binary recursive relations. In case of non-degenerate sequences we give bounds for the highest index for which a term can be 0 and bounds on the growth order of the absolute values of the terms, both only in terms of the two initial values, which is a novel feature. Some of these bounds are best possible apart from a multiplicative constant.
\end{abstract}

\maketitle

\section{Introduction}

In this paper we study integers represented  by sequences defined as follows. Let $A, B, P$ and $Q$ be integers with $AB\neq 0$ and one of $P,Q$ being non-zero.
Let $(u_n)_{n=0}^\infty$ be given by $u_0=P$, $u_1=Q$ and
\begin{equation} \label{rec}
u_n=Au_{n-1}- Bu_{n-2}\ \ \ (n\geq 2).
\end{equation}
We call such a sequence a BRIG sequence (Binary Recursive Integrally Generated sequence). Properties of such sequences have been studied in \cite{bhpr}, \cite{gl}, \cite{kov}, \cite{pet}, \cite{stta}, \cite{yh}. Computational aspects can be found in section `Binary recursive sequences' of the website of SAGE \cite{sage}.

A BRIG sequence is called a Lucas sequence of first kind if $u_0 = 0$, $u_1 =1$, a Lucas sequence of second kind if $u_0 = 2$, $u_1 = A$. Lucas sequences are well studied in the literature (see e.g. \cite{bw} and the references given there). They include the sequences of the Fibonacci numbers for $(A,B) = (1,-1)$, of the Pell numbers for $(A,B)=(2,-1)$, of the Jacobsthal numbers for $(A,B)=(1,-2)$ and of the Mersenne numbers for $(A,B)=(3,2)$.

Write $f(x)=x^2-Ax+B$ for the characteristic polynomial of $u$. Let $\alpha,\beta$ be its roots. Throughout the paper we shall assume that $|\alpha|\geq |\beta|$. Assuming $\alpha \neq \beta$ we have
\begin{equation} \label{exprun}
u_n = a \alpha^n - b \beta^n \ \ \ (n\geq 0),
\end{equation}
where
\begin{equation} \label{ab}
a = \frac{Q - P \beta} {D}, \ \ b = \frac{Q - P \alpha}{D}, \ \ D= |\alpha - \beta |= \sqrt{A^2-4B}.
\end{equation} 
In particular, the Lucas sequences of first and second kind are given by 
$$U_n = \frac {\alpha^n - \beta^n}{\alpha - \beta}, \ \ \  V_n = \alpha^n + \beta^n,$$
respectively.
If $\alpha$ is real, we say that we are in the real case, otherwise in the non-real case. We use $c_1, c_2, \ldots $ for effectively computable absolute constants $>1$. Unless stated otherwise we assume that $ \alpha / \beta$ is not a root of unity. This together with the requirements $abAB \neq 0$ and one of $P,Q$ is non-zero is expressed by writing that the sequence is non-degenerate.

In Section \ref{mres} we formulate our principal results, which concern zeroes and growth properties of BRIG sequences. Similar results are known in the literature: for zeroes of linear recurrence sequences see e.g. \cite{av} and the references there, for growth properties of such sequences see \cite {st76} (cf. \cite{ss} and \cite{htl}). In particular, it follows from Theorem \ref{thm21} of \cite{htl} that $|U_n| \geq \frac 12 |\alpha|^{n-2}$ in the real case and $|U_n| \geq |\alpha|^{n-c_1 (\log n)^2}$ in the non-real case for $n \geq 5$, independent of the chosen non-degenerate Lucas sequence. We use this result in Section \ref{zerbrig} to show that for any non-degenerate BRIG sequence there is at most one index $n$ such that $u_n=0$, that $u_n = 0$ implies that $n < 10 \log |Q| + c_2$, and that this bound is best possible apart from a multiplicative factor. In Sections \ref{lbbrig} and \ref{lbnr} we derive lower bounds for the terms of BRIG-sequences comparable to those for Lucas sequences in \cite{htl}. The lower bound for $|u_n|$, of the form $u_n > c_3^n$, holds in the real case for $n$ greater than an explicit expression only depending on $P$ and $Q$, but in the non-real case on an expression depending on $B, P$ and $Q$. 
The proofs in the real case are elementary, those in the non-real case are based on estimates of linear forms in logarithms.

The novelty of our results are that they are explicit, in some cases are quantitatively sharp, and are qualitatively best possible: they depend only on those parameters which cannot be avoided.


\section{Main results} \label {mres}

In this section we use the following notation and assumptions. Let $A, B, P$ and $Q$ be integers with $abAB \neq 0$. Define the sequence $(u_n)_{n=0}^{\infty}$ by \eqref{rec}. Write $x^2-Ax+B = (x- \alpha)(x- \beta)$ with $|\alpha| \geq |\beta| >0$. We assume that $\alpha / \beta$ is not a root of unity. Then $u_n$ is given by \eqref{exprun} for all $n$ with $a,b$ as in \eqref{ab}; we also define $D$ as in \eqref{ab}. Numbers $c_1, c_2, \ldots $ denote effectively computable absolute constants $>1$.

For simplicity, we shall also assume that $PQ\neq 0$. Clearly, if $PQ=0$ then $(u_n)_{n=m}^\infty$ is a multiple of a Lucas sequence for $m=0$ or $1$, and the results follow easily from those given in \cite{htl}. Beside this, we discuss the excluded sequences in Remark 3.3 at the end of Section 3.

Theorem \ref{thm22} concerns the zero terms in BRIG sequences. The zero-multiplicity (i.e. the number of zero terms) of a recurrence sequence has been studied in many papers, see e.g. \cite{av} and the references there. We believe that part a. follows easily from results and observations made in the literature, but we could not find a concrete reference. So, as it is important for parts b. and c., we formulate and later prove this statement as well.

\begin{theorem} \label{thm22} {\rm a.} For given $A,B,P,Q$ there is at most one integer $k$ such that $u_k=0.$\\
{\rm b.} If $u_k = 0$ and $A^2 > 4B$, then $k< 9 \log|Q| + 12$. \\
{\rm c.} If $u_k=0$ and $A^2<4B$, then $k < 10 \log \max(|Q|,2)$ for $k > c_4$.
\end{theorem}

\begin{theorem} \label{thm23}
Suppose $A^2 > 4B$.\\
If $|A-D| \geq 6\left|\frac QP\right|$, then for $n \geq 6 |\frac QP | + 6$ we have $$|u_n| \geq |Q|\left( \frac{|\alpha|} {2}\right)^{n-2} \ \ {\rm and} \ \ |u_n| \geq |Q|\left( \frac {\sqrt{5}}{2} \right)^n.$$
If $|A-D| < 6|\frac QP|$, then  for $n \geq (18 + 7 \log |Q|) \max \left(1, \left|\frac {Q}{P}\right|\right)$ we have $$|u_n|  \geq  \min\left(\frac{1}{5|P|}, \frac{1} {22 |Q|}\right) |\alpha|^{n-2} \ \ {\rm and} \ \ |u_n| \geq \min\left(\frac{1}{14|P|},\frac {1}{36|Q|}\right)\left(\frac{1+ \sqrt{5}}{2}\right)^n .$$
\end{theorem}
\noindent 
The proof provides sharper bounds.

\begin{theorem} \label{thm24} Suppose $A^2<4B$.
For $$n> c_5 \log(B|P|+|Q|) (\log \log (B|P|+|Q|))^2,$$ we have
$$|u_n| \geq |\alpha|^{2n/3}  > 1.25^n.$$
\end{theorem}

We conjecture that also in the non-real case it should be possible to give a lower bound on $n$ which depends only on $P$ and $Q$.


\section{On the zeros of BRIG sequences} \label{zerbrig}
We use the notation of the first paragraph of Section \ref{mres} and apply Lemma \ref{thm21} to prove Theorem \ref{thm22}.

\begin{lemma} [\cite{htl}, Theorem 2.1]
\label{thm21}
Define the sequence $(U_n)_{n=0}^{\infty}$ by $U_0 = 0, U_1=1$ and $U_n = AU_{n-1} - BU_{n-2}$ for $n=2,3, \dots $. Suppose $(U_n)$ is non-degenerate and $n \geq 2$. Then
\begin{equation} \label{Ugt}
|U_n| \geq \frac 12 |\alpha|^{n-2}~~~{\rm if}~~~B<0,
\end{equation}\begin{equation} \label{UBneg}
|U_n| \geq  |\alpha|^{n-1}~~~{\rm if}~~~0 < 4B <A^2,
\end{equation}
and 
\begin{equation} \label{Ust}
 |U_n| \geq |\alpha|^{n-c_1 (\log n)^2} ~~~{\rm if}~~~A^2 < 4B.
\end{equation}
\end{lemma}
Remark 1 of \cite{htl} states that $|\alpha| \geq 2$ in the real case except for $(A,B) = (1,-1)$, and that $|\alpha| \geq \sqrt{2}$ in the non-real case.

\vskip.2cm
\noindent{\it Proof of Theorem} \ref{thm22}. To prove part a., assume that $k$ is the smallest index with $u_k=0$. Then $u_{k+1}\neq 0$, since otherwise
we would have $P=Q=0$. Observe that then the sequence $(u_n)_{n=k}^\infty$ is $u_{k+1}$ times a Lucas sequence, and we have
$$
u_n=u_{k+1}\frac{\alpha^{n-k}-\beta^{n-k}}{\alpha-\beta}\ \ \ (n\geq k).
$$
Hence the claim follows.

Now we turn to the proof of parts b. and c.
In the sequel we assume that $A>0$ as we may by replacing $A$ by $-A$ if $A<0$. Observe that then if $\alpha$ and $\beta$ are real, we have $\alpha>|\beta|>0$.

Let $d$ be the largest integer such that both $d~|~A$ and $d^2 |~B$. We claim that we may assume that $d=1$ without loss of generality.
Indeed, we have $d^{n-1}~|~u_n$ by induction. Put $A' = \frac Ad, B'= \frac{B}{d^2}$ and
$$ P' = u'_0 = dP,~~~Q' = u'_1 = Q, ~~~u'_n = A'u'_{n-1} - B'u'_{n-2} \ \ (n=2,3,\dots).$$
Then $P'$ and $Q'$ are integers and $$u'_n = \frac{u_n}{d^{n-1}} \in \mathbb{Z} \ \ {\rm for~all} ~n.$$
Moreover $u_k=0$ if and only if $u_k'=0$.

We may further assume that $\gcd(P,Q)=1$. Indeed, otherwise we can divide every term of the sequence by $\gcd(P,Q)$ to get smaller integer values of $|P|$ and $|Q|$. 
In the sequel of this section we assume $d=1, \gcd(P,Q) = 1$, $u_k=0$ for some integer $k>1$.

By induction it is obvious that, for given $A,B$ and arbitrary $n \geq 2$, $u_n$ is a homogeneous linear combination of $P$ and $Q$ with coefficients $c_P(n,A,B)$ and $c_Q(n,A,B)$ depending only on $n, A$ and $B$. Observe that
$$
c_P(0,A,B)=c_Q(1,A,B)=1,\ \ \ c_P(1,A,B)=c_Q(0,A,B)=1,
$$
$$
c_P(n+2,A,B)=Ac_P(n+1,A,B)-Bc_P(n,A,B)\ \ \ (n\geq 0),
$$
$$
c_Q(n+2,A,B)=Ac_Q(n+1,A,B)-Bc_Q(n,A,B)\ \ \ (n\geq 0),
$$
and as one can readily check,
$$
c_P(n+1,A,B)=-Bc_P(n,A,B)\ \ \ (n\geq 1).
$$
From this we easily get that $c_P(n,A,B)\neq 0$, $c_Q(n,A,B)\neq 0$ for $n\geq 2$.
Thus $\frac {P}{Q}$ is uniquely determined by $k, A$ and $B$. We fix $k, A$ and $B$ and shall compute the corresponding fraction $\frac PQ$. 

We define the Lucas sequence $(U_n)$ by $U_0=0, U_1=1$ and $U_n = AU_{n-1}-BU_{n-2}$ for $n=2,3, \dots $. Further we define a sequence $(u''_n)_{n=0}^{k}$ by $u''_n = B^n U_{k-n}$. Then $ u''_k = B^k U_0 = 0$. Observe that 
$$u''_n -Au''_{n-1} + Bu''_{n-2} = B^n U_{k-n} -AB^{n-1}U_{k-n+1} + B^{n-1} U_{k-n+2}=$$ 
$$= B^{n-1}(U_{k-n+2} - AU_{k-n+1} + BU_{k-n} )= 0$$
for $n=2,3, \dots, k$. It follows, by induction for $n=k, k-1, k-2, \dots, 0$, that $u_n = \frac {u_{k-1}}{u''_{k-1}} u''_n$ for $n=0,1, \dots, k$. Thus the sequence $(u_n)_{n=0}^k $ is a constant multiple of the sequence $(u''_n)_{n=0}^k$ and 
\begin{equation}\label{expPQ}
\frac PQ = \frac {u_0}{u_1} = \frac {u''_0}{u''_1} = \frac {U_k}{BU_{k-1}}.
\end{equation}
\vskip.1cm
\noindent {\bf Example 3.1.} Let $A=3, B=6, k=5.$ We obtain $U_0=0, U_1=1, U_2 = 3, U_3 = 3, U_4 = -9, U_5 = -45$. Hence $u''_0 = -45, u''_1 = -54, u''_2 = 108, u''_3 = 648, u''_4 =1296, u''_5=0$. 
The sequence $(u''_n)$ satisfies the recurrence $u''_n = Au''_{n-1} - Bu''_{n-2}$ for all $n \geq 2$, and, apart from a common multiplicative factor for all terms, this is the unique sequence for these values of $A, B$ and $k$ with $u_5=0$. Here $\frac PQ = \frac {-45}{-54} = \frac 56.$ 
\vskip.1cm

Let $g:= \gcd(A,B).$ Then, by $d=1$, for every prime $p$ which divides $g$ we have $v_p(A) \geq v_p(B)=1$, hence $v_p(g)=1$. By induction on $n$ we obtain $v_p(U_{2n+1}) = n$ and $v_p(U_{2n}) \geq n$ for all $n$. Therefore $$v_p(\gcd(U_{2n},U_{2n+1})) = v_p(\gcd(U_{2n+1},U_{2n+2})) = n$$ for all $n$. If $p\nmid B$, $p\mid U_n$, $p\mid U_{n-1}$, then, by induction using \eqref{rec}, $p\mid U_1=1$, a contradiction. If $p\mid B$, $p\mid U_n$ for any $n\geq 2$, then $p\mid A$ and therefore $p\mid\mid g$. Thus 
\begin{equation} \label{eqD}
\gcd(U_{2n},U_{2n+1}) = \gcd(U_{2n+1},U_{2n+2}) = \prod_{p | g} p^n = g^n
\end{equation} for all $n$.
\vskip.1cm
\noindent {\bf Example 3.2.} Let $A=15, B=10.$ Then $d=1, g=5$ and $$(U_n) = \{0, 1, 15, 215, 3075, 43975, 628875, ...\}.$$ For the sequence $(\gcd(U_{n},U_{n+1}))$ we find $1, 1, 5, 5, 25, 25, ...$ . \qed
\vskip.3cm


We assume $u_k=0$ with $k>6$ without loss of generality.
Note that, for all $A$ and $B$,
\begin{equation} \label{boundQ}
|Q| =  \frac {|BU_{k-1}|} { \gcd(U_k,BU_{k-1})} \geq \frac {|U_{k-1}|} {g^{\frac{k-1}{2}}}.
\end{equation}

If $A^2>4B>0$, then $A \geq 3$, $\alpha \geq \frac{A+1}{2} \geq 2$ and $g \leq A$. Thus, by \eqref{boundQ} and \eqref{UBneg}, if $\alpha \geq \frac 52$,
$$
|Q|  \geq \frac {|U_{k-1}|} {g^{\frac{k-1}{2}}} \geq \frac {\alpha^{k-2}}{(2 \alpha)^{\frac {k-1}{2}}} \geq \frac 12 \left(\frac{\alpha}{2}\right)^{\frac{k-3}{2}} \geq \frac 12 \left(\frac{5}{4}\right)^\frac{k-3}{2},
$$
and if $2 \leq \alpha <\frac 52$,
$$
|Q| \geq \frac {|U_{k-1}|} {g^{\frac{k-1}{2}}} \geq \frac {\alpha^{k-2}}{(2\alpha -1)^{\frac{k-1}{2}}} \geq \frac {\alpha ^{k-2}}{\left(\frac 85 \alpha\right)^{\frac{k-1}{2}}} = \left(\frac 58\right)^{\frac{k-1}{2}}\alpha^{\frac{k-3}{2}} \geq \frac 58 \left(\frac 54\right)^{\frac{k-3}{2}}.
$$

If $B<0$, then $g \leq A <\alpha < D$. 
If $A \geq 2$, then $\alpha \geq 1 + \sqrt{2}$ and, by \eqref{boundQ} and \eqref{Ugt},
$$ |Q| \geq \frac{U_{k-1}}{g^{\frac{k-1}{2}}} \geq \frac{ \frac12 \alpha^{k-3}}{\alpha^{\frac {k-1} {2}}}
= \frac12  \alpha^{\frac{k-5}{2}} \geq \frac12 (1+\sqrt{2})^{\frac{k-5}{2}}.$$
If $A=1$, then $g=1$ and $\alpha \geq \frac 12 (1+\sqrt{5})$. It follows that
$$|Q| \geq U_{k-1} \geq \frac{\alpha^{k-3}}{2}  \geq \frac12 \left(\frac{1+\sqrt{5}}{2}\right)^{k-3}.$$
Thus, if $A^2>4B$, then $|Q| \geq \frac 12 \min\left((1.25)^{\frac{k-3}{2}}, (1+\sqrt{2})^{\frac {k-5}{2}}\right)$, hence
\begin{equation} \label{pos}
k < 9 \log |2Q| + 5 \leq 9 \log |Q| + 12.
\end{equation}

If $A^2<4B$ and $B \leq 5$, then the only pairs $(A,B)$ with $\gcd(A,B)>1$ are $(2,2), (3,3)$ and $(2,4)$. These three pairs lead to degenerate cases.
In all other cases $A$ and $B$ are coprime. The case $B=1$ leads to a degenerate case. Thus $g=1$ and $|\alpha| =  \sqrt{B} \geq \sqrt{2}$. So we obtain for $B \leq 5$, by \eqref{Ust}, 
\begin{equation} \label{les5}
|Q| = B|U_{k-1}| \geq \ |\alpha|^{2+ k-1- c_6 (\log k)^2} \geq 2^{\frac k2 - c_6(\log k)^2}. 
\end{equation}

If $A^2<4B$ and $B>5$, then we have, by \eqref{Ust} and $g \leq |A| \leq 2 \sqrt{B} = 2 |\alpha|$,
$$
|Q| = \frac {B|U_{k-1}|} { \gcd(U_k,BU_{k-1})} \geq \frac {|\alpha|^{k+1 - c_6 (\log k)^2}} {g^{\frac k2-1}} \geq \frac {|\alpha|^{k - c_6(\log k)^2}} {(2|\alpha|)^{\frac k2-1}} 
\geq \frac{|\alpha|^{\frac k2-c_6 (\log k)^2}}{2^{\frac k2 -1}}.
$$
Therefore, since $|\alpha| = \sqrt{B} \geq \sqrt{6}$,
\begin{equation} \label{asmb}
|Q| \geq \frac{(\sqrt{6})^{\frac k2 - c_6 (\log k)^2}}{2^{\frac k2}} \geq (\sqrt{1.5})^{\frac k2-c_7 (\log k)^2}.
\end{equation}
This yields $k < 10 \log |Q|$ for $k > c_8$.

The combination of inequalities \eqref{pos}, \eqref{les5} and \eqref{asmb} completes the proof of Theorem \ref{thm22}.
\qed

\vskip.2cm
\noindent {\bf Remark 3.1.} Let $(u_n)$ be a non-degenerate BRIG sequence with $PQ\neq 0$ and $u_k=0$. Then $P = u_0 \neq 0$. Define $v_n = P^{n-1}u_n$ for all $n \geq 0$. Then $(v_n)_{n=0}^{\infty}$ is an integer sequence with $v_0=1, v_1 = Q$ and $v_n =APv_{n-1} -BP^2v_{n-2}$ for $n \geq 2$ which satisfies $v_k=0$ if and only if $u_k=0$. So if we apply the bound on $k$ to the sequence $(v_n)$ instead of the sequence $(u_n)$, there is no dependence on $P$.
\vskip.2cm
\noindent {\bf Remark 3.2.} 
The following example shows that in the non-degenerate case a logarithmic upper bound in $Q$ for $k$ is best possible.
Let $k$ be an arbitrary integer $>2$.
Set $$P=u_0=2^k-1,~ Q= u_1=2^k-2,~ A =3,~ B = 2.$$
Then $u_n = 2^k - 2^{n}$ for all $n$ and $u_k=0$ so that $k = \frac{\log (Q+2)}{ \log 2} > 1.44 \log Q$.
\vskip.2cm
\noindent{\bf Remark 3.3.}
We have assumed that $abABPQ\neq 0$ and that $\frac {\alpha}{\beta}$ is not a root of unity. Here we check which cases in Theorem \ref{thm22} were excluded by these restrictions.

As mentioned before, the case $PQ=0$ immediately reduces to the Lucas sequences treated in \cite{htl}. So we may suppose that both are not zero.
If the sequence $(u_n)$ is of the form $a \alpha^n$ for $n>0$, then it is constant $0$ or it is non-zero for $n>0$. This happens if $B=0$, hence $\alpha = 0$ or $\beta =0$, and if $ab=0$. From here on we assume $abB \neq 0$.
If $A=0$, then $\alpha = - \beta \not= 0$ and $\frac{\alpha}{\beta} = -1$ is a root of unity.  

If $\frac{\alpha }{ \beta}$ is a root of unity, then we distinguish between $\alpha = \beta$ and $\alpha \neq \beta$. If $\alpha = \beta$, that is $A^2=4B$, then $u_n = nQ(\frac A2)^{n-1} - (n-1)P(\frac A2)^n~~ {\rm for~all}~n.$ Hence $u_k= 0$ if and only if $2kQ=(k-1)PA$. Thus, for given $A,P,Q$,  there is at most one $k>1$ with $u_k=0$. 

Suppose $\alpha / \beta$ is a root of unity with $\alpha \neq\beta$. Then $(\alpha / \beta)^m = 1$ for some $m \leq 6$. 
Therefore the sequence is periodic with period $m$ apart from a multiplicative factor $a\alpha^n$. Thus either there are infinitely many zeros or none. E.g. the choice $A=B=P=Q=1$ leads to a sequence $(u_n)$ with period 6 and $u_n=0$ if and only if $n \equiv 0 \pmod 3$.
\vskip.1cm
For given $A,B$ we give values for $P=P_k,Q=Q_k$ for sequences with $u_k = 0$:\\
For $u_2=0$ we can choose $P_2=A, Q_2=B$.\\
For $u_3 = 0$ we can choose $P_3=A^2-B, Q_3=AB$.\\
In general we can choose $P_{m+1} = AP_m - Q_m, Q_{m+1} = BP_m$ for $m \geq 3$. This follows by induction in view of $AQ_{m+1}-BP_{m+1} = BQ_m$. The new sequence 
$(u_n^{(m+1)})_{n=1}^{\infty}$ equals $B$ times the old sequence $(u_n^{(m)})_{n=0}^{\infty}$ and therefore $u_{m+1}^{(m+1)}=0$.

\section{The growth of BRIG-sequences in the real case} \label{lbbrig}
In this section we use the notation as in the first paragraph of Section \ref{mres} and prove Theorem \ref{thm23}. 
Cam Stewart (\cite{st76} p. 33, cf. \cite{ss} Lemma 5) proved the following result:

\begin{theorem}
\label{sthm}
Suppose $\alpha /\beta$ is not a root of unity.
Then there exist computable numbers $C_{1}$ and $C_{2}$ depending only on $a$ and $b$ such that 
$$|u_n| \geq |\alpha|^{n - C_{1} \log n}~~~~(n \geq C_{2}).$$
\end{theorem}
\noindent In this section we shall prove a similar result with $C_{1}$ and $C_{2}$ depending only on $P$ and $Q$. Since we may divide all the terms by the same integer, we may assume without loss of generality that $\gcd(P,Q)=1$, $P \geq 0$, and that if $P=0$ then $Q=1$. As before, we assume $A>0$. Hence, as now $\alpha$ and $\beta$ are real, we have $\alpha>|\beta|>0$.

By \eqref{exprun} and \eqref{ab} we have
\begin{equation} \label{abc}
 \alpha = \frac{A + D}{2}, \beta =  \frac{A - D}{2}, a = \frac{Q - P\beta }{D}, b = \frac{Q - P\alpha}{D}.
 \end{equation}
 \noindent
If $D \in \mathbb{Z}$, then $b/a$ satisfies the linear equation $$(PA-PD-2Q)x -  (PA+PD-2Q) =0.$$ 
If $D \notin \mathbb{Z}$, then $b/a$ is a root of the irreducible polynomial \\
\begin{equation} \label{Dab}
D^2(bx-a)(ax-b)= 
\end{equation}
$$
\left(Q^2-PQA+P^2\frac{A^2-D^2}{4}\right)(x^2 + 1)- \left(2Q^2-2PQA+P^2\frac{A^2+D^2}{2}\right) x.
$$
This polynomial has integer coefficients, since $A^2-D^2=4B$ is divisible by $4$ and therefore $A^2+D^2$ is even. In both cases we find, denoting the canonical height of the integer $r$ by $H(r)$,
\begin{equation} \label{c1c2}
H \left(\frac{b}{a} \right) \leq 2Q^2 + 2|PQ|A +P^2~\frac{A^2+D^2}{2} \leq 2\left(|Q|+ |P|~\frac{A+|D|}{2}\right)^2-1.
\end{equation}

We shall use that, for any algebraic number $\gamma$,
\begin{equation} \label{height}
\frac{1}{H(\gamma) + 1} < |\gamma| < H(\gamma)+1
\end{equation}
see \cite{cij}, Lemma 1.2. For the convenience of the reader we give the short proof. If ${\frac{1}{2}}<\gamma<2$ then the claim is trivial. We may assume $|\gamma| \geq 2.$ Let $T(z) = t_dx^d + t_{d-1}x^{d-1} + \ldots + t_1x+t_0$ be the minimal primitive polynomial of $\gamma$ over $\mathbb{Z}$ and $h$ its height. Then
$$ |\gamma|^d \leq |t_d \gamma^d| \leq h(|\gamma|^{d-1} + \ldots + |\gamma| +1) < h |\gamma|^{d-1} (1 - |\gamma|^{-1})^{-1}$$
This implies $|\gamma|<h+1$. The minimal polynomial of $\frac{1}{\gamma}$ is the reciprocal polynomial, hence $\frac{1}{|\gamma|} < h+1$. So \eqref{height} follows.

Furthermore we shall use that
\begin{equation} \label{basic}
|a \alpha^n - b \beta^n | \ \geq | a \alpha^n |~ \left( 1 - \left| \frac {b} {a} \left(\frac {\beta}{\alpha}\right)^n\right| \right)
\end{equation}
and that, by \eqref{abc},
\begin{equation} \label {casea1}
| a | =  \frac {|Q - P \beta |} {|D| }  \cdot \frac {| Q - P \alpha |} {| Q - P \alpha |} = \frac {|Q^2-APQ+BP^2|}{|D|~|Q - P \alpha |} \geq \frac {2}{|D|( 2|Q| + |P|(|A+D|))}.
\end{equation}

 \vskip.2cm
 \noindent {\it Proof of Theorem} \ref{thm23}.
Suppose $A^2 >4B$ and $A>0$. This implies that $\alpha$ and $\beta$ are real numbers with $\alpha > |\beta |$.  If $B>0$, then $B \geq 1$, $A \geq 3, D\geq 1, \alpha \geq 2$. If $B<0$, then either $A=1, B=-1$, $D = \sqrt{5}$, $\alpha = \frac 12 (1+ \sqrt{5})$ or $D \geq \sqrt{8}$, $\alpha \geq 2$. Note that $n\geq 7$ in view of $|\frac QP| > 0.$

We distinguish the following four cases: \\(a) $ A-D \geq 6|\frac QP|$, (b) $ D-A \geq 6|\frac QP|$, (c) $|A-D| <  6|\frac QP|$ and $ \\A+D \geq 9|\frac QP|$, (d) $|A-D| < 6|\frac QP|$ and $ A+D < 9|\frac QP|$.
\vskip.3cm
\noindent {\bf Case (a)}: $ A-D \geq 6|\frac QP|$. We have $A\geq 3, B \geq 1, D \geq 1$ and
$$u_n = a \alpha^n - b \beta^n = \frac {Q - P\beta }{D} \left(\frac {A+D}{2}\right)^n -  \frac {Q - P\alpha }{D} \left(\frac {A-D}{2}\right)^n=$$
$$ \frac{1}{2^{n+1}D} \{(2Q-PA+PD)(A+D)^n - (2Q-PA-PD)(A-D)^n\}=$$
$$-\frac{1}{2^{n}} \{(PA-PD-2Q) \left(\frac{(A+D)^n - (A-D)^n}{2D}\right) - P(A-D)^n\}.$$
\vskip.2cm
We use that, for $x>y>0$,
\begin{equation} \label{devxy}
(x+y)^n - (x-y)^n = 2 \sum_{k=1,~ k~{\rm odd}}^n {n \choose k} x^{n-k}y^k \geq 2nx^{n-1}y.
\end{equation} 
We obtain, for $n \geq 7,$
$$2^{n} |u_n| \geq (|P|(A-D)-2|Q|) \cdot nA^{n-1} - |P| (A-D)^n \geq$$
$$|P|A^{n-1}(A-D)\left(n -\frac n3 -1\right) \geq  (4n-6) |Q| A^{n-1} \geq 22|Q|A^{n-1}.$$ 
Thus, for $n \geq 7$, by $A \geq 3, A > \alpha$,
\begin{equation} \label{casa}
|u_n| \geq 11|Q| \left(\frac A2\right)^{n-1} \geq 11|Q|\left(\frac {\alpha}{2}\right)^{n-1}~~{\rm and}~~ |u_n| \geq 7|Q|\left(\frac 32\right)^{n}.
\end{equation}

\vskip.1cm
\noindent {\bf Case (b)}: $ D-A \geq 6|\frac QP|$. Then $A \geq 1, B \leq -1, D\geq \sqrt{5}$ and, as in Case (a),
\begin{equation} \label{expun}
u_n=  \frac{1}{2^{n+1}D}\{(PD-PA+2Q)(D+A)^n + (PD+PA-2Q)(A-D)^{n} \}.
\end{equation}
We distinguish between two subcases.

\noindent (b1) If $n$ is even, then 
$$|u_n| \geq  \frac{|P|}{2^{n+1}D} \left\{\left(D-A-2\left|\frac QP\right|\right)(D+A)^n + \left(D+A-2\left|\frac QP\right|\right)(D-A)^n\right\} \geq$$
$$\frac{|P|}{2^{n+1}D} \left\{4\frac {|Q|}{|P|}(D+A)^n  \right\} > \frac {|Q|(D+A)^{n-1}}{2^{n-1}}.$$
Thus, if $n$ is even, then
\begin{equation} \label{casb1}
|u_n| \geq  |Q| \alpha^{n-1} ~~{\rm and} ~~|u_n| \geq 0.6 |Q|\left(\frac {1+\sqrt{5}}{2}\right)^n.  
\end{equation}

\noindent (b2) Now suppose $n$ is odd and $D-A \geq 6|\frac QP|, n \geq 6 | \frac QP | +3$, hence $n \geq 7$. Therefore, by \eqref{expun} and \eqref{devxy},
$$ \frac{2^{n+1}D}{|P|} |u_n| \geq \left(D-A -2|\frac QP|\right)\{(D+A)^n - (D-A)^n\} - \left(2A+ 4|\frac {Q}{P}|\right) (D-A)^{n} \geq$$
$$\frac 23 (D-A)\cdot 2nD^{n-1}A - \left(2A+4 |\frac QP|\right)(D-A)D^{n-1}  \geq  $$
$$ \left(\frac {4n}{3} -2 -4|\frac QP|\right) A(D-A)D^{n-1} \geq \frac {2n}{3} A \cdot 6 |\frac QP | D^{n-1}.$$
Thus we obtain, by $n \geq 7$, $D>\alpha$, $D \geq \sqrt{5}$,
\begin{equation} \label{casb2}
|u_n| \geq   \frac {1}{2} nA|Q|(\frac D2)^{n-2}>  \frac 72 A|Q|(\frac {\alpha}{2})^{n-2}~~ {\rm and }~~ |u_n| \geq 2.8|Q|\left(\frac{\sqrt{5}}{2}\right)^n.
\end{equation}
\noindent 
Combining formulas \eqref{casa}, \eqref{casb1}, \eqref{casb2} and using $n \geq 7$, we obtain the first statement of Theorem \ref{thm23}.

\vskip.3cm
\noindent {\bf Case (c)}: $|D-A| < 6|\frac QP|, A+D \geq 9 |\frac QP|$. Then, by $D>0$, \eqref{c1c2} and \eqref{height}, 
$$\left| \frac{b}{a} ~ \left(\frac{\beta}{\alpha}\right)^n\right| \leq 2\left(|Q|+|P|\frac{A+D}{2}\right)^2 \cdot \frac{36Q^2} {(A+D)^2P^2} \cdot \left( \frac{|A-D|}{A+D}\right)^{n-2} \leq $$
$$18 \left( \frac{2Q^2}{|P|(A+D)} + |Q|\right)^2 \left(\frac {2}{3}\right)^{n-2} \leq 27~Q^2\left(\frac {2}{3}\right)^{n-2}< \frac 12$$
if $n > 2 + \frac {\log 54 + 2\log |Q|}{ \log 1.5}$. If $n$ is that large, then, by \eqref{casea1} and \eqref{basic},
$$|u_n| \geq \frac{1}{D(2|Q|+|P|(A+D))} \left(\frac{A+D}{2}\right)^n  \geq \frac{9}{44|P|} \left(\frac{A+D}{2}\right)^{n-2}.$$
Thus, if $n > 12 + 5 \log |Q|$, then
\begin{equation} \label{casc}
|u_n| \geq \frac{1}{5|P|}\alpha^{n-2}  \geq \frac{1}{14|P|}  \left( \frac {1 + \sqrt{5}}{2} \right)^n.
\end{equation}

\vskip.1cm
\noindent {\bf Case (d)}: $|D-A| < 6|\frac QP|, A+D < 9|\frac QP|$. 
Hence $1 \leq A\leq 7.5|\frac QP|, \\1 \leq D \leq 7.5 |\frac QP|$ and, by \eqref{casea1},
\begin{equation} \label{inequn}
2^{n} |u_n| \geq  \frac{2(A+D)^n}{D(2|Q|+|P|(A+D))} - \frac {2|Q|+|P|(A+D)}{2D} |A-D|^n \geq 
\end{equation}
$$\frac {(A+D)^n} {5.5|Q|D} - 5.5\frac{|Q|}{D} |A-D|^n \geq \frac{(A+D)^n}{5.5|Q|D} \left(1 - 31Q^2 \left( \frac{|A-D|}{A+D} \right)^n\right).$$
Since $|D-A| < 6|\frac QP|$, $A \geq 1$, $D \geq 1$, we have
\begin{equation} \label{ADPQ}
\frac {A+D}{|A-D|} \geq 1 + \frac {2}{|A-D|}  > 1 + \frac {|P|}{|3Q|}.
\end{equation}
Suppose $|P| \geq |Q|$. Then 
\begin{equation} \label{PgrQ}
\log \left(1+ \frac{|P|}{3|Q|}\right) \geq \log \frac 43 \geq \frac 27.
\end{equation}
Suppose, on the contrary, $|P|<|Q|$. Then, by the monotonicity of $ \frac {\log(1+x)}{x}$ for $x>0$,
\begin{equation} \label{PsmQ}
\log \left( 1 + \frac {|P|}{3|Q|}\right) \geq \frac{|P|}{3|Q}| \cdot 3 \log \frac43 > \frac 27 \cdot \frac{|P|}{|Q|}.
\end{equation}
On combining \eqref{ADPQ}, \eqref{PgrQ} and \eqref{PsmQ} we obtain that for $$n \geq (18 + 7 \log |Q|) \max\left(1, \left|\frac QP \right|\right)$$ we have 
$$n \log \frac{A+D}{|A-D|} > \frac 27 n \min\left(1, \frac{|P|}{|Q|}\right) \geq 5 + 2\log |Q| \geq \log 62 + 2 \log |Q|.$$
Thus $31Q^2 (\frac{|A-D|}{A+D})^n < \frac 12$ and, by \eqref{inequn},
\begin{equation} \label{casd}
|u_n| \geq \frac {1} {11|Q|D} \left(\frac{A+D}{2}\right)^n \geq \frac {1}{22|Q|} \alpha^{n-1} \geq \frac{1}{36|Q|} ~\left(\frac {1+\sqrt{5}}{2}\right)^n.
\end{equation}

\noindent Combining formulas \eqref{casc} and \eqref{casd} we obtain the second statement of Theorem \ref{thm23}.\qed

\vskip.1cm
\noindent {\bf Remark 4.1.} One may expect that the lower bound for $|u_n|$ is linear in $Q$. This is true for the bound in case $|D-A| \geq 6|Q/P|$, but not in the other case. The reason is technical. In the former case $A$ or $D$ is dominating and this results in linearity in $Q$. If $D-A < 6|Q/P|$, then the terms $a\alpha^n$ and $b \beta^n$ may have comparable size and we apply the rather rough estimates \eqref{c1c2} and \eqref{casea1} to distinguish between them.

\section{The growth of BRIG-sequences in the non-real case} \label{lbnr}
In this section we use the notation of the first paragraph of Section \ref{mres} and prove Theorem \ref{thm24}. 

Suppose $A^2 < 4B$. In this case $\alpha$ and $\beta$ are complex conjugates and $B>0$. If $B=1$, then 
$\alpha / \beta$ is a root of unity. 
We assume that $\alpha / \beta$ is not a root of unity and therefore $B \geq 2$, $|\alpha| \geq \sqrt{2}$.
We recall that $(\alpha - \beta)^2 = A^2 - 4B$ and that $\frac{\beta}{\alpha}$ is a zero of the quadratic polynomial with integer coefficients
$B( x - \frac{\alpha}{\beta}) ( x - \frac{ \beta}{\alpha}) = Bx^2 - (A^2-2B)x + B.$
Since $-2B \leq A^2-2B \leq 2B$, it follows that the height of $\frac{\beta}{\alpha}$ is at most $2B$. 
The first inequality of \eqref{c1c2} implies that
$$ H \left(\frac{b}{a} \right) \leq 2BP^2 +2A|PQ| +2Q^2 \leq 2(|P|\sqrt{B} + |Q|)^2.$$

We use a theorem of Baker \cite{ba77} in the following form:
\begin{lemma}  [\cite{spts}, Theorem A] \label{BS}
Let $\alpha_1, \dots, \alpha_n$ be non-zero algebraic numbers. Let $K$ be their splitting field over $\mathbb{Q}$. Put $E = [K:Q]$. We denote by $A_1, \dots, A_n$ upper bounds for the heights of $\alpha_1, \dots, \alpha_n$, respectively,  where we assume that $A_j \geq 2$ for $1 \leq j \leq n$. Write 
$$\Omega' = \prod_{j=1}^{n-1} \log A_j, ~~~~\Omega = \Omega' \log A_n.$$
There exist computable numbers $c_9$ and $c_{10}$ such that the inequalities
\begin{equation} \label{baker77}
0 < \mid \alpha_1^{b_1} \cdots \alpha_m^{b_m} - 1 \mid < \exp( -(c_9mE)^{c_{10}m} \Omega \log \Omega' \log M)
\end{equation}
have no solution in integers $b_1, \dots, b_n$ with absolute values at most \\$M (\geq 2)$.
\end{lemma}

\noindent{\bf Remark 5.1.} A result of Matveev \cite{mat} provides a sharper estimate. However, the above statement suffices for our purposes.

\vskip.2cm

\noindent We apply Lemma \ref{BS} with $m=E=2,  \alpha_1 = \beta / \alpha, \alpha_2 = b/a, b_1=n,\\ b_2 =1$ and obtain, when $u_n \neq 0$,
\begin{equation} \label{compl} 
\left|\frac{b}{a} \left(\frac{\beta}{\alpha}\right)^n -1\right| \geq \exp(-c_{11} \log(|P| \sqrt{B} + |Q|) \log(2B) \log \log(2B) \log n).
\end{equation}
We have $$|\alpha| = \sqrt{|\alpha \beta|} = \sqrt{B}~~ {\rm and}~~ a = \frac {2Q-PA}{2\sqrt{A^2-4B}} + \frac P2.$$
The first term of $a$ is purely imaginary, the second is real, so that $|a| \geq |\frac P2| \geq \frac 12.$
Combining this with \eqref{compl} we get, by \eqref{basic}, 
$$
|u_n| \geq \frac 12 B^{n/2} \exp\left(- c_{12} \log(B|P|+|Q|) \log(2B) \log \log(2B) \log n\right). 
$$
Thus there is an absolute constant $c_{13}$ such that if 
\begin{equation} \label{n0}
n > n_0 := c_{13} \log(B|P|+|Q|) (\log \log (B|P|+|Q|))^2
\end{equation}
we have
$$
|u_n| \geq B^{n/3} = |\alpha|^{2n/3} \geq 2^{n/3}.
$$

\end{document}